# Ensemble transform Kalman-Bucy filters

Javier Amezcua [1 #]
jamezcua@atmos.umd.edu

Kayo Ide [1 2 3 4]
ide@umd.edu

Eugenia Kalnay [1 2 3]
ekalnay@atmos.umd.edu

Sebastian Reich [5]
sereich@rz.uni-potsdam.de

[1] Department of Atmospheric and Oceanic Science. University of Maryland, College Park, MD, USA.
[2] Earth System Science Interdisciplinary Center. University of Maryland, College Park, MD, USA.
[3] Institute for Physical Science and Technology. University of Maryland, College Park, MD, USA.
[4] Center for Scientific Computation and Mathematical Modeling. University of Maryland, College Park, MD, USA.
[5] Insitut für Mathematik. Am Neuen Palais 10, D-14469. Universität Potsdam, Potsdam, Germany.
[#] The principal author's current affiliation is: Department of Meteorology, University of Reading, UK.

# Abstract

Two recent works have adapted the Kalman-Bucy filter into an ensemble setting. In the first formulation, BGR09, the ensemble of perturbations is updated by the solution of an ordinary differential equation (ODE) in pseudo-time, while the mean is updated as in the standard KF. In the second formulation, BR10, the full ensemble is updated in the analysis step as the solution of single set of ODEs in pseudo-time. Neither requires matrix inversions except for the frequently diagonal observation error covariance.

We analyze the behavior of the ODEs involved in these formulations. We demonstrate that they stiffen for large magnitudes of the ratio of background to observational error covariance, and that using the integration scheme proposed in both BGR09 and BR10 can lead to failure. An integration scheme that is both stable and is not computationally expensive is proposed. We develop transform-based alternatives for these Bucy-type approaches so that the integrations are computed in ensemble space where the variables are weights (of dimension equal to the ensemble size) rather than model variables.

Finally, the performance of our ensemble transform Kalman-Bucy implementations is evaluated using three models: the 3-variable Lorenz 1963 model, the 40-variable Lorenz 1996 model, and a medium complexity atmospheric general circulation model (AGCM) known as SPEEDY. The results from all three models are encouraging and warrant further exploration of these assimilation techniques.

# 1. Introduction

Two recently proposed methods implement the (continuous) Kalman-Bucy Filter (KBF) (Kalman and Bucy, 1961) in an ensemble framework. Both formulate the analysis step through an ordinary differential equation (ODE) representation in "pseudo-time". In BGR09 (Bergemann *et al.*, 2009) only the ensemble of perturbations is updated in this way; in BR10 (Bergemann and Reich, 2010a) the full ensemble is updated as the solution of a single set of ODEs. One major advantage of BR10 is that it allows "mollification" of the observational increments (Bergemann and Reich, 2010b), thus maintaining the balance of the analysis state at least as well as the widely used Intermittent Analysis Update (Bloom *et al.*, 1996). Furthermore, both formulations are suitable for extensions that deal with non-Gaussian uni-modal and multi-modal ensemble distributions (Reich, 2011), and they are a natural option in the case when observations are assimilated in an almost continuous fashion (Bergemman and Reich, 2011).

An essential challenge in both BGR09 and BR10 is to integrate their respective ODEs in pseudo-time in an efficient and affordable manner. We perform a detailed analysis of the behavior of these ODEs. In particular we evaluate the robustness of the Bucy-type formulations in the presence of large perturbation growth that results from infrequent observations, which is usually accompanied by non-linear behavior (Kalnay *et al.*, 2007a; Yang and Kalnay, 2009a). It is shown that in this case the ODEs involved in the Bucy-type formulations stiffen; to handle this problem a diagonally semi-implicit (DSI) modification of the Euler Forward scheme is presented. Moreover, we propose transform-based formulations for both BGR09 and BR10; under these schemes the calculations are performed in the ensemble space rather than in the state space.

Finally, the transform-based Bucy-type formulations are tested with three models of increasing complexity. This paper is organized as follows: Section 2 provides background information, with 2.1 briefly reviewing the KF, KBF and EnKF and 2.2 giving a summary of BGR09 and BR10. Section 3 contains the analysis of the ODEs and proposes an improved numerical integration method for the Bucy-type methods. In section 4 we develop our transform-based alternatives. In section 5 we perform experiments with the Lorenz 1963 model (L63, section 5.1), the Lorenz 1996 model (L96, section 5.2) and an AGCM of medium complexity (SPEEDY model, section 5.3). For all our experiments, we use as a benchmark for comparison the local ensemble transform Kalman filter (LETKF; Hunt *et al.*, 2007). A summary of the results and conclusions is presented in section 6.

## 2. Background

### 2.1. Kalman filter and Kalman-Bucy filter

The KF (Kalman, 1960) deals with the filtering process of a time-discrete linear dynamical system:

$$\mathbf{x}_t = \mathbf{A}\mathbf{x}_{t-1} \qquad \mathbf{x} \in \Re^N \tag{1}$$

$$\mathbf{y}_t = \mathbf{H}\mathbf{x}_t + \mathbf{v}_t \qquad \mathbf{y} \in \Re^L \tag{2}$$

Equation (1) represents the linear evolution of the state variables $\mathbf{x} \in \Re^N$ by action of the state transition matrix $\mathbf{A} \in \Re^{N \times N}$, assuming a perfect model. Equation (2) defines the observations $\mathbf{y} \in \Re^L$ as a linear combination of state variables transformed by the matrix $\mathbf{H} \in \Re^{L \times N}$ plus a stochastic noise term $\mathbf{v}_t \sim N(\mathbf{0}, \mathbf{R})$ that represents the observational error, where $\mathbf{R} \in \Re^{L \times L}$ is the observational error covariance (usually assumed to be diagonal in atmospheric applications). The KF solves this estimation problem as a two-step process. The first is the forecast of both the model state $\hat{\mathbf{x}}$ and its error covariance matrix $\mathbf{P} \in \Re^{N \times N}$:

$$\mathbf{P}_t^b = \mathbf{A}\mathbf{P}_{t-1}^a \mathbf{A}^{\mathrm{T}} \tag{3}$$

$$\hat{\mathbf{x}}_t^b = \mathbf{A}\hat{\mathbf{x}}_{t-1}^a \tag{4}$$

The notation $\hat{\mathbf{x}}$ indicates this is an estimator of the truth $\mathbf{x}$. The superscript $b$ stands for background (forecast) before the assimilation of the observations at that time. The superscript $a$ stands for analysis, where this information has been included. The second step is the analysis and assimilates the observation at the time in the background; hereafter we drop the time index $t$ unless otherwise noted:

$$\mathbf{P}^a = (\mathbf{I} - \mathbf{KH})\mathbf{P}^b = \left((\mathbf{P}^b)^{-1} + \mathbf{H}^T \mathbf{R}^{-1} \mathbf{H}\right)^{-1} \tag{5}$$

$$\hat{\mathbf{x}}^a = \hat{\mathbf{x}}^b - \mathbf{K}(\mathbf{H}\hat{\mathbf{x}}^b - \mathbf{y}) \tag{6}$$

The matrix $\mathbf{K} \in \mathfrak{R}^{N \times L}$ is known as the Kalman gain, and is given by:

$$\mathbf{K} = \mathbf{P}^b \mathbf{H}^T \left(\mathbf{H} \mathbf{P}^b \mathbf{H}^T + \mathbf{R}\right)^{-1} = \mathbf{P}^a \mathbf{H}^T \mathbf{R}^{-1} \tag{7}$$

The time-continuous KBF (Kalman and Bucy, 1961) deals with the filtering process of the time-continuous linear dynamical system:

$$\frac{d\mathbf{x}}{dt} = \mathbf{F}\mathbf{x} \tag{8}$$

where $\mathbf{F}$ is the linear dynamical model. Time-continuous observations $\mathbf{y}$ of the time-continuous model state $\mathbf{x}$ have the same form as in the KF. In this case, the forecast and the analysis steps are merged together into a system of ODEs (see for example Simon, 2006):

$$\frac{d\mathbf{P}}{dt} = \mathbf{F}\mathbf{P} + \mathbf{P}\mathbf{F}^T - \mathbf{P}\mathbf{H}^T \mathbf{R}^{-1} \mathbf{H}\mathbf{P} \tag{9}$$

$$\frac{d\hat{\mathbf{x}}}{dt} = \mathbf{F}\hat{\mathbf{x}} - \mathbf{P}\mathbf{H}^T \mathbf{R}^{-1} (\mathbf{H}\hat{\mathbf{x}} - \mathbf{y}) \tag{10}$$

For a time-discrete filtering of the system given by (1) and (2), the KF analysis (5)-(7) can be reformulated using the ODE approach similar to the time-continuous KBF given by (9) and (10). We call this filtering the time-discrete KBF. The forecast step of the time-discrete KBF follows the time-discrete KF, i.e., (3) and (4). If the dynamics is governed by the time-continuous model as in (8) while the observations occur at discrete times, then $\mathbf{A}$ in (4) is given by the resolvent of (8). The reformulated analysis step at a fixed (and discrete) time $t$ is given by the time-continuous ODEs in *pseudo-time s* spanning over $0 \leq s \leq 1$:

$$\frac{d\mathbf{P}}{ds} = -\mathbf{P}\mathbf{H}^T \mathbf{R}^{-1} \mathbf{H}\mathbf{P} \tag{11}$$

$$\frac{d\hat{\mathbf{x}}}{ds} = -\mathbf{P}\mathbf{H}^T \mathbf{R}^{-1}(\mathbf{H}\hat{\mathbf{x}} - \mathbf{y}) \tag{12}$$

with the initial conditions $\mathbf{P}(0) = \mathbf{P}^b$ and $\hat{\mathbf{x}}(0) = \hat{\mathbf{x}}^b$. At $s$=1, the solution of the ODEs coincide with the analysis, i.e., $\mathbf{P}^a = \mathbf{P}(1)$ and $\hat{\mathbf{x}}^a = \hat{\mathbf{x}}(1)$. That is, starting from the background at the beginning of the pseudo-time window, (11) and (12) give the analysis at the end (Bergemann *et al.*, 2009). The time-discrete KBF only involves the (usually simple) matrix inversion $\mathbf{R}^{-1}$.

## 2.2. Ensemble Kalman Bucy filters (EnKBFs): BGR09 and BR10

Many geophysical data assimilation systems of interests make observations at discrete times. The remaining of this paper focuses on the time-discrete filtering. Dynamics of these systems are often nonlinear:

$$\mathbf{x}_t = f(\mathbf{x}_{t-1}) \tag{13}$$

where *f* represents the nonlinear model operator. To effectively deal with the nonlinearity of the dynamics for high-dimensional systems while taking advantage of the KF analysis to assimilate the observations, the ensemble Kalman Filter (EnKF) uses an ensemble of finite size $M$ as follows. Using the matrix representation of the ensemble of the model state:

$$\bar{\bar{\mathbf{X}}} = [\mathbf{x}_1 \mid \mathbf{x}_2 \mid \cdots \mid \mathbf{x}_M] \qquad \bar{\bar{\mathbf{X}}} \in \Re^{N \times M} \qquad (14)$$

where the subscripts denote the ensemble members, the sample mean can be written as:

$$\bar{\mathbf{x}} = M^{-1} \bar{\bar{\mathbf{X}}} \mathbf{1} \qquad (15)$$

using $\mathbf{1} \in \Re^M$ as a column vector of ones. The estimator $\hat{\mathbf{x}}$ that appeared in the KF analysis equations is replaced by this sample mean $\bar{\mathbf{x}}$. An ensemble of perturbations can be obtained by subtracting the sample mean column from each ensemble member:

$$\mathbf{X} = [\mathbf{x}_1 - \bar{\mathbf{x}} \mid \mathbf{x}_2 - \bar{\mathbf{x}} \mid \cdots \mid \mathbf{x}_M - \bar{\mathbf{x}}] = \bar{\bar{\mathbf{X}}}(\mathbf{I} - \mathbf{U}) \qquad (16)$$

where $\mathbf{U} = M^{-1} \mathbf{1}\mathbf{1}^T$, $\mathbf{U} \in \Re^{M \times M}$. Then the sample covariance can be obtained by:

$$\mathbf{P} = (M-1)^{-1} \mathbf{X}\mathbf{X}^T = (M-1)^{-1} \bar{\bar{\mathbf{X}}}(\mathbf{I} - \mathbf{U})\bar{\bar{\mathbf{X}}}^T \qquad (17)$$

The forecast step of the EnKF forwards $\mathbf{X}^a$ at time *t-1* to $\mathbf{X}^b$ at time *t*. This step is straightforward, given by the ensemble forecast of (13). The analysis step obtains $\mathbf{X}^a$ from $\mathbf{X}^b$ at time *t* by assimilating observations $\mathbf{y}$ in a way that the covariance (17) and the mean (15) of the analysis ensemble are consistent with the KF counterparts in (5) and (6). This step is not unique, and can be done either in a stochastic fashion ('perturbed observations' method: Burgers *et al.*, 1998; Houtekamer and Mitchell, 1998) or in a deterministic fashion (Tippett *et al.*, 2003).

Both BGR09 and BR10 belong to the family of the deterministic ensemble square root filters. They differ from others by adapting the time-discrete KBF (11) and (12) into the analysis step. We call them ensemble KBFs (EnKBFs). As we will see in section 4, our new extensions of BGR09 and BR10 are based on the local ensemble transform Kalman filter (LETKF; Hunt *et al.*, 2007), which is another ensemble square root filter. The analysis step of BGR09 and BR10 is formulated as follows:

Starting from (11), the first part of the analysis step in BGR09 is to update the ensemble perturbations from $\mathbf{X}^b$ to $\mathbf{X}^a$ by solving:

$$\frac{d\mathbf{X}}{ds} = -\frac{1}{2(M-1)}\mathbf{X}\mathbf{X}^T\mathbf{H}^T\mathbf{R}^{-1}\mathbf{H}\mathbf{X} \tag{18}$$

in pseudo-time s over $0 \leq s \leq 1$ with initial condition $\mathbf{X}(0) = \mathbf{X}^b$. At $s=1$, we obtain the analysis ensemble of perturbations $\mathbf{X}^a = \mathbf{X}(1)$. Once $\mathbf{X}^a$ is obtained, the analysis mean $\bar{\mathbf{x}}^a$ can be computed using (6) with $\mathbf{K} = (M-1)^{-1}\mathbf{X}^a\mathbf{X}^{aT}\mathbf{H}^T\mathbf{R}^{-1}$. The full analysis ensemble $\bar{\bar{\mathbf{X}}}^a$ is constructed from $\mathbf{X}^a$ and $\bar{\mathbf{x}}^a$.

In contrast, BR10 combines (11) and (12) into a single set of ODEs and directly obtains the analysis ensemble $\mathbf{X}^a$ by solving:

$$\frac{d\bar{\bar{\mathbf{X}}}}{ds} = -\frac{1}{2(M-1)}\mathbf{X}\mathbf{X}^T\mathbf{H}^T\mathbf{R}^{-1}\left[\mathbf{H}\bar{\bar{\mathbf{X}}}(\mathbf{I}+\mathbf{U}) - 2\mathbf{y}\mathbf{1}^T\right] \tag{19}$$

in pseudo-time $s$ over $0 \leq s \leq 1$ with the initial condition $\bar{\bar{\mathbf{X}}}(0) = \bar{\bar{\mathbf{X}}}^b$. At $s=1$, we obtain the analysis snsemble $\bar{\bar{\mathbf{X}}}^a = \bar{\bar{\mathbf{X}}}(1)$. A derivation of (19) can be found in the appendix. In BR10, this appears as the equivalent form of the gradient of a cost function.

A fundamental difference between BGR09 and BR10 is that BGR09 computes the analysis mean $\bar{\mathbf{x}}^a$ separately using (6). Accordingly, observations $\mathbf{y}$ themselves impact only the analysis mean $\bar{\mathbf{x}}^a$ but not the analysis ensemble perturbation $\mathbf{X}^a$. In BGR10, observations $\mathbf{y}$ affect both $\mathbf{X}^a$ and $\bar{\mathbf{x}}^a$ because of the nonlinearity in (19). This difference can be important for numerical implementations as the update of the ensemble mean will be affected by the chosen time-stepping. It is worth mentioning that neither BGR09 nor BR10 requires any matrix inversion except $\mathbf{R}^{-1}$, which is frequently diagonal and usually can be pre-computed.

## 3. Numerical integration in the EnKBF

To develop efficient time-stepping approximations for the EnKBFs, we need to study the qualitative solution behavior of the ODE (18) in BGR09. If accuracy alone is of concern, explicit time-stepping approaches are the methods of choice. In fact, BGR09 used the Euler forward method with 4 steps and found that adding more steps did not significantly improve the root mean square error (RMSE) of the analysis. However, in this section, we demonstrate that the ODEs in (11) can become stiff depending on the ratio of $\mathbf{P}$ and $\mathbf{R}$. Two main controlling factors of this ratio are the frequency of observations and length of the assimilation window. When the ratio is large, the forward Euler approximation of (11) loses stability unless very small time-steps are applied. We address this challenge with a simple diagonally semi-implicit modification to the Euler method.

Consider the analytical solution of the Bucy equation for the covariance in pseudo-time (11) (see Simon, 2006):

$$\mathbf{P}(s) = \mathbf{P}^b \left( \mathbf{H}^T \mathbf{R}^{-1} \mathbf{H} \mathbf{P}^b s + \mathbf{I} \right)^{-1} \tag{20}$$

For illustration purposes, let $\mathbf{H} = \mathbf{I}$ and consider the scalar case of (20), i.e. $N = 1$. It is convenient to write it as the following fraction, which is the ratio of the covariance at any moment $s$ in pseudo-time and the background covariance:

$$\frac{P_1(s)}{P_1^b} = \frac{1}{\beta s + 1} \tag{21}$$

In (21), $\beta = P_1^b / \sigma_1^2$ is the ratio of background covariance to observational error covariance (variance, in this scalar example). The behavior of (21) depends upon the magnitude of this ratio as illustrated in figure 1. For small values of $\beta$ (upper-right curves), the function behaves close to linear, since $P_1(s)/P_1^b \approx 1 - \beta s$ when $\beta < 1$. It is possible to get an accurate numerical approximation of these curves with the Euler forward method and using few steps. Nonetheless, for large values of $\beta > 1$ (lower curves), the problem becomes increasingly stiff and an explicit method such as Euler forward is no longer effective without significantly increasing the number of steps. Returning to the general case with several variables and an arbitrary $\mathbf{H}$, the expression for $\beta$ becomes:

$$\beta = \left| \mathbf{Y}^{bT} \mathbf{R}^{-1} \mathbf{Y}^b \right| / (M - 1) \tag{22}$$

where $\mathbf{Y}^b = \mathbf{H}\mathbf{X}^b$ is the mapping of the ensemble of perturbations into the observation space and $||\ ||$ denotes the spectral matrix norm, i.e., the maximum singular value. $\mathbf{P}^b$ and therefore $\beta$ depend on the length of the forecast window, as well as on the sensitivity of the observation operator $\mathbf{H}$. For short enough windows, $\mathbf{P}^b$ remains smaller than the observational error covariance $\mathbf{R}$ (hence $\beta < 1$), while for long windows it can become considerably larger (hence $\beta > 1$). As an illustration, we consider an assimilation experiment using LETKF and the Lorenz 1963 model for both short ($\beta_s$) and long ($\beta_l$) windows (defined by having observations every 8 or 25 steps respectively; the details are explained in section 5.1). The empirical cumulative distribution functions (ECDF) of $\beta_s$ and $\beta_l$ are shown in figure 2. From this figure we see that the ECDF for $\beta_s$ is an order of magnitude smaller than the one for $\beta_l$. For short windows $\beta_s < 10^{-1}$ occurs about 45% of the cycles, and $\beta_s > 1$ only 15% (in fact $\max(\beta_s) = 4.81$). By contrast, for long windows $\beta_l < 10^{-1}$ only 1% of the times, $\beta_l > 1$ for about 60% of the cases, and $\max(\beta_l) = 187.07$. Hence, for the latter case the ODEs involved in the Bucy-type formulations are bound to stiffen.

To tackle this stiffness problem, we use a diagonally semi-implicit (DSI) Euler approximation in place for the Euler forward time scheme used in BR09:

$$\mathbf{X}_{k+1} = \mathbf{X}_k - \frac{\Delta s}{2}\mathbf{P}_k \mathbf{H}^T (\Delta s \mathbf{D}_k + \mathbf{R})^{-1}\mathbf{H}\mathbf{X}_k \qquad (23)$$

where $\mathbf{P}_k = (M-1)^{-1}\mathbf{X}_k \mathbf{X}_k^T$, $\mathbf{D}_k$ is the diagonal of $\mathbf{H}\mathbf{P}_k\mathbf{H}^T$, and $k$ is the pseudo-time step index, The inversion introduced is performed on a diagonal matrix and adds a negligible cost compared to the Euler Forward method. In fact, replacing $\mathbf{R}^{-1}$ by $(\Delta s \mathbf{D}_k + \mathbf{R})^{-1}$ does not change the order of accuracy of the Euler forward approximation (Hairer and Wanner, 1991). It does improve, however, the stability and hence it is useful when the step size $\Delta s$ is restricted by the method's stability rather than the accuracy. The DSI method of (25) falls into the category of Rosenbrock methods with inexact Jacobian; it is also called W methods in Hairer and Wanner (1991). Moreover, this method can handle non-diagonal $\mathbf{R}$ by replacing $(\Delta s\, \mathbf{D}_k + \mathbf{R})^{-1}$ with $\left(diag\left(\mathbf{HP}_k\mathbf{H}^T\mathbf{R}^{-1}\Delta s + \mathbf{I}\right)\right)^{-1}\mathbf{R}^{-1}$.

The resulting update for the ensemble mean is:

$$\bar{\mathbf{x}}_{k+1} = \bar{\mathbf{x}}_k - \Delta s \mathbf{P}_k \mathbf{H}^T \left(\Delta s\, \mathbf{D}_k + \mathbf{R}\right)^{-1}\left[\mathbf{H}\bar{\mathbf{x}}_k - \mathbf{y}\right] \qquad (24)$$

which becomes the standard Kalman update over a time-interval $\Delta s$ if $\mathbf{D}_k$ is replaced by the full matrix $\mathbf{HP}_k\mathbf{H}^T$. Then the DSI scheme of BR10 given by (19) is:

$$\bar{\bar{\mathbf{X}}}_{k+1} = \bar{\bar{\mathbf{X}}}_k - \frac{\Delta s}{2} \mathbf{P}_k \mathbf{H}^T \left(\Delta s\, \mathbf{D}_k + \mathbf{R}\right)^{-1}\left[\mathbf{H}\bar{\bar{\mathbf{X}}}_k(\mathbf{I}+\mathbf{U}) - 2\mathbf{y}\mathbf{1}^T\right] \qquad (25)$$

To complement the DSI method, we choose a sequence of pseudo-time steps with increasing size that sum to one. For example a sequence of 7 steps will be {1/16,1/16,1/8,1/4,1/4,1/4}, a sequence of 8 steps will be {1/32,1/32,1/16,1/8,1/4,1/4,1/4} and so forth. This sequence was found heuristically to work well; the rationale behind it is that the fastest change of equations (18) and (19) occurs at the beginning of pseudo-time, just around $s = 1/\beta$, so shorter steps are required there. The advantages of using this integration scheme with respect to the original Euler forward method will be illustrated in section 5. A further investigation on the selection of the step size could be developed in the future.

**4. Ensemble Transform Kalman-Bucy filters**

In the ensemble-transform approach, the EnKF expresses the analysis ensemble of perturbations as a product of the background ensemble of perturbations $\mathbf{X} \in \Re^{N \times M}$ in the model space and a matrix of weights in the ensemble space $\mathbf{W} \in \Re^{M \times M}$, i.e. $\mathbf{X}^a = \mathbf{X}^b \mathbf{W}^a$. This transforms the analysis step in the model-state space to the ensemble space, which has a much lower dimension. This approach leads to the ensemble transform KF (ETKF; Bishop *et al*, 2001; Wang *et al*, 2004; Ott *et al*, 2004) and local ETKF (LETKF: Hunt *et al*, 2007). By incorporating the ensemble-transform approach, we develop two EnKF methods that are stemmed from BGR09 and BR10 and refer them as the ensemble transform KBFs (ETKBFs).

Having the weights available in the ensemble space can be convenient for some existing applications, e.g., accurate low-resolution analyses by weight interpolation (Yang *et al.*, 2009a), a no-cost smoother (Kalnay *et al*, 2007b), forecast sensitivity to observations without adjoint model (Liu and Kalnay, 2008, Li et *al.*, 2009), and Running in Place/Quasi Outer-Loop (Kalnay and Yang, 2010; Yang and Kalnay, 2009). These techniques rely on the fact that the weights at the analysis time are valid through the entire forecast window (Kalnay *et al*, 2007b, Yang and Kalnay, 2009)and been shown to work effectively for data sparse situations that give rise to nonlinear perturbations (Kalnay and Yang, 2010; Yang *et al*, 2009).

Application of the ensemble-transform approach to the EnKBFs is straightforward. Starting from BGR09 (18) and using $\mathbf{X}(s) = \mathbf{X}^b \mathbf{W}(s)$, we find an ODE for the evolution of matrix of weights $\mathbf{W}(s)$ in pseudo-time. This is the basis for the ensemble transform Kalman-Bucy filter (ETKBF):

$$\frac{d\mathbf{W}}{ds} = -\frac{1}{2(M-1)} \mathbf{W}\mathbf{W}^T \mathbf{Y}^{bT} \mathbf{R}^{-1} \mathbf{Y}^b \mathbf{W} \tag{26}$$

Solving this in pseudo-time *s* over $0 \leq s \leq 1$ with the initial condition $\mathbf{W}(0) = \mathbf{I}$, we obtain the analysis weight matrix $\mathbf{W}^a = \mathbf{W}(1)$ at *s*=1. Hence, if $\mathbf{W}$ satisfies (26), then $\mathbf{X}$ satisfies (18). It can be shown that the solution for an infinite number of steps in pseudo-time is equivalent to the form used in LETKF (Appendix B). This equivalence no longer holds when the observation operator is not linear. In the Bucy-type formulations, the forward operator could be adjusted/linearized after each pseudo-time step; the study of this case is beyond the scope of this work but constitutes an interesting are of research.

The corresponding DSI integration scheme for the ETKBF is:

$$\mathbf{W}_{k+1} = \mathbf{W}_k - \frac{\Delta s}{2} \tilde{\mathbf{P}}_k^T \mathbf{Y}^{bT} \left(\Delta s \tilde{\mathbf{D}}_k + \mathbf{R}\right)^{-1} \mathbf{Y}^b \mathbf{W}_k \tag{27}$$

where $\tilde{\mathbf{P}}_k = (M-1)^{-1} \mathbf{W}_k \mathbf{W}_k^T$ and $\tilde{\mathbf{D}}_k = diag\left(\mathbf{Y}^b \tilde{\mathbf{P}}_k \mathbf{Y}^{bT}\right)$.

Similarly BR10 (19) leads to the Direct Ensemble Transform Kalman-Bucy filter (DETKBF) by letting $\overline{\overline{\mathbf{W}}} \in \mathfrak{R}^{M \times M}$ transform the background ensemble into the analysis ensemble, i.e. $\overline{\overline{\mathbf{X}}}^a = \overline{\overline{\mathbf{X}}}^b \overline{\overline{\mathbf{W}}}^a$. We note that the full ensemble space matrix $\overline{\overline{\mathbf{W}}}$ is different from the perturbation matrix $\mathbf{W}$ of the ETKBF. The DETKBF obtains the matrix of weights $\overline{\overline{\mathbf{W}}}$ by solving:

$$\frac{d\overline{\overline{\mathbf{W}}}}{ds} = -\frac{1}{2(M-1)} \mathbf{W}\mathbf{W}^T \overline{\overline{\mathbf{Y}}}^{bT} \mathbf{R}^{-1} \left[\overline{\overline{\mathbf{Y}}}^b \overline{\overline{\mathbf{W}}}(\mathbf{I} + \mathbf{U}) - 2\mathbf{y}\mathbf{1}^T\right] \tag{28}$$

in pseudo-time s over $0 \leq s \leq 1$ with the initial condition $\overline{\overline{\mathbf{W}}}(0) = \mathbf{I}$ where $\overline{\overline{\mathbf{Y}}}^b = \mathbf{H}\overline{\overline{\mathbf{X}}}^b$ is the mapping of the (full) background ensemble into observations space and $\mathbf{W} = \overline{\overline{\mathbf{W}}}(\mathbf{I} - \mathbf{U})$. At *s*=1, we obtain the analysis weight matrix $\overline{\overline{\mathbf{W}}}^a = \overline{\overline{\mathbf{W}}}(1)$. Hence, if $\overline{\overline{\mathbf{W}}}(1)$ is the solution of (28), then $\overline{\overline{\mathbf{X}}}^a = \overline{\overline{\mathbf{X}}}^b \overline{\overline{\mathbf{W}}}^a$ is the solution to (25). The corresponding DSI integration scheme is:

$$\overline{\overline{\mathbf{W}}}_{k+1} = \overline{\overline{\mathbf{W}}}_k - \frac{\Delta s}{2} \tilde{\mathbf{P}}_k \overline{\overline{\mathbf{Y}}}^{bT} \left( \Delta s \, \tilde{\mathbf{D}}_k + \mathbf{R} \right)^{-1} \left[ \overline{\overline{\mathbf{Y}}}^b \overline{\overline{\mathbf{W}}}_k (\mathbf{I} + \mathbf{U}) - 2\mathbf{y}\mathbf{1}^T \right] \tag{29}$$

where $\tilde{\mathbf{D}}_k$ is defined as in (27).

## 5. Experiments with three models

We conduct three sets of experiments. Each one illustrates a particular aspect: the value of the DSI method (section 5.1), the use of localization and adaptive inflation (section 5.2) and practical applicability to a realistic model (section 5.3). For the three models we use identical twin experiments. Most practical data assimilation systems need two basic algorithms, localization and inflation, to attain a successful performance. When the ensemble size is much smaller than not only the dimension of the model state ($M \ll N$) but also the number of the positive Lyapunov exponents, straightforward application of any EnKF may lead to unreliable correlation estimations especially at long distance. The gridpoint **R**-localization is a simple yet powerful technique to handle this challenge for the EnKFs with the ensemble-transform approach; in this scheme, an independent analysis is carried out for every single grid point using observations within a certain distance and assuming that the observation error increases with the distance to the grid point (see Hunt *et al.*, 2007 and Greybush *et al.*, 2011, for details). Underestimation of the background ensemble perturbation may also occur due to small $M$ and other sources of model imperfection. For the EnKFs with the ensemble-transform approach, an adaptive inflation scheme (Miyoshi, 2011) addresses this issue. These techniques are employed in the experiments in sections 5.2 and 5.3.

## 5.1. The Lorenz 1963 model

Our first set of experiments relies on L63, a nonlinear 3-variable model widely used in evaluating data assimilation schemes because of its challenging properties near regime changes (e.g. Miller *et al.*, 1994; Evensen, 1998; Evans *et al.*, 2004). The system of nonlinear coupled ODEs describing its evolution is:

$$\begin{aligned}\dot{x}^{(1)} &= p(x^{(2)} - x^{(1)}) \\ \dot{x}^{(2)} &= x^{(1)}(r - x^{(3)}) - x^{(2)} \\ \dot{x}^{(3)} &= x^{(1)}x^{(2)} - bx^{(3)}\end{aligned} \qquad (30)$$

The standard values are used for the parameters: $p=10$, $r=28$ and $b=8/3$, which result in a chaotic behavior with two regimes in a very well-known butterfly-shaped fractal attractor in the phase space. The model is integrated with the Runge-Kutta 4$^{th}$ order method using a time step of $\Delta t = 0.01$.

We use similar settings similar to those of Kalnay *et al.* 2007a and Miller *et al.* 1994. The "observations" are generated by adding a random noise $N(\mathbf{0}, \mathbf{R} = 2\mathbf{I})$ to the nature run. Two cases are considered: "frequent" observations at every 8 time steps and "infrequent" observations at every 25 time steps. Frequent observations lead to assimilation windows in which the perturbations grow essentially linearly; this roughly corresponds to a 6-hr assimilation cycle in an atmospheric global circulation model. With infrequent observations, the perturbations grow nonlinearly and their distribution is no longer Gaussian so that Kalman filtering becomes less accurate (Kalnay *et al*, 2007a). For both assimilation frequencies a sample of $10^6$ analysis cycles was used for the results reported.

For assimilation, we use $M = 3$, the smallest possible size in this model. The ensemble members are initialized by adding random noise to the truth with the same covariance as the observational error.

Multiplicative covariance inflation (Anderson, 2001) is used. We choose $\delta > 0$ such that the background ensemble of perturbations is multiplied by a factor $\mathbf{X}^b \to \mathbf{X}^b(1+\delta)$, equivalent to multiplying the background covariance matrix by $\mathbf{P}^b \to \mathbf{P}^b(1+\delta)^2$. For both observation frequencies, we varied $\delta$ in order to minimize the analysis RMSE. For the Bucy-type formulations we also search for the minimum number of pseudo-time steps that lead to an analysis RMSE comparable to that of the LETKF.

For the frequent observations case, the values of $\delta$ are taken from $\{0.01, 0.02, \ldots, 0.1\}$ (in Kalnay *et al.*, 2007a, the optimal value was found to be $\delta = 0.04$). For ETKBF and DETKBF it is found that using less than 3 steps for the integration leads to a poor performance. The filters start performing well with 3 steps, and with just 5 steps the performances of both Bucy-type filters are the same as that of LETKF. The computing time for the three methods is indistinguishable. In the left panel of figure 3 we depict the performance of the formulations integrated using 5 steps. The three filters show similar behavior with respect to the inflation parameter. From $\delta = 0$ to about $\delta \sim 0.03$, the performance of the filters improves fast as inflation increases. After this value, there is an optimal performance region for the three filters; the lowest analysis RMSE values are almost identical: for ETKF 0.3108 ($\delta = 0.07$), for ETKBF 0.3064 ($\delta = 0.06$) and for DETKBF 0.3163 ($\delta = 0.06$). Beyond this region, the covariance inflation becomes excessive and the filter begins to lose skill slowly.

For the infrequent observations case, inflation values are taken from $\{0.1, 0.2, \ldots, 0.9\}$ (in Kalnay *et al.* 2007a, the optimal inflation was $\delta = 0.39$). Our first experiments use Euler forward to integrate both Bucy-type formulations. As expected from the stiffening, this scheme fails with a number of steps of $O(10)$. A large number of pseudo-time steps (~70 for ETKBF and ~300 for DETKBF) are necessary to achieve a performance similar to LETKF, and occasional failure is still observed. Therefore we switch to the DSI integration method. If we use uniform pseudo-time steps, we find that at least 30 steps are needed for ETKBF to achieve the performance of the LETKF (for DETKBF this number was ~50). Switching to the variable time stepping discussed in section 3 reduces these numbers to 8. The results are depicted in the right panel of figure 3. Again, for the 3 formulations a rapid reduction in RMSE is observed as one increases the value of inflation before $\delta \sim 0.03$; then an optimal inflation region is found. The lowest analysis RMSE values are: for ETKF 0.7544 ($\delta = 0.4$), for ETKBF 0.7664 ($\delta = 0.5$) and for DETKBF 0.7612 ($\delta = 0.5$), about 1% larger for both Bucy-based formulations than for LETKF (but with comparable computing time).

**5.2. The Lorenz 1996 model**

To test the effect of localization, the 40-variable L96 model (Lorenz 1996; Lorenz and Emanuel, 1998). This periodic model is described by the set of differential equations:

$$\dot{x}^{(q)} = \left(x^{(q+1)} - x^{(q-2)}\right)x^{(q-1)} - x^{(q)} + F \tag{31}$$

for $q = 1,\ldots,40$, with $x^{(0)} = x^{(40)}$, $x^{(-1)} = x^{(39)}$ and $x^{(41)} = x^{(1)}$. $F = 8$ is a forcing term. The attractor of this model has a fractal dimension of about 27 and 13 positive Lyapunov exponents (Lorenz, 2005). It does not have regime transitions as L63.

The model is integrated with a Runge-Kutta 4$^{th}$ order method and a time step of $\Delta t = 0.025$ units. Observations are taken every 2 time steps which is roughly equivalent to a 6 hours window in an atmospheric general circulation model (Lorenz and Emmanuel, 1998). We observe every other grid point with an observational error covariance $\mathbf{R} = \mathbf{I}$ as in BR10.

Two ensemble sizes are considered: $M = 10$ and $M = 15$. R-localization is used with a well-known compact support function: equation (4.10) in Gaspari and Cohn 1999 with $c = \lambda/\sqrt{.3}$, with a localization radius of $\lambda = 4$. Values $\lambda = \{1.5, 2, 2.5, \ldots, 9\}$ are tested and this value was experimentally found to minimize the RMSE (Amezcua, 2012). For multiplicative inflation we avoid manual tuning of $\delta$ by using adaptive multiplicative covariance inflation (Miyoshi, 2011). This scheme -which uses the diagnostic relationships of Desroziers *et al*, 2005- is tailored for R-localization and estimates a time-evolving individual inflation parameter each gridpoint.

By varying the number of pseudo-time steps ETKBF and DETKBF, the performance with 3 steps was comparable to that of the LETKF and after 4 steps we found only marginal improvements. Table 1 shows the results of this experiment for the three methods (columns) and the two ensemble sizes (rows); A sample of $10^6$ assimilation cycles is used to compute the average analysis RMSE, average ensemble spread and average inflation parameter (averaged also over the 40 gridpoints). The numbers in parenthesis correspond to the standard deviations of the reported parameters. The overall performance is indistinguishable among the 3 methods, and the computational time is the same.

The DSI method is effective when initializing the background ensemble without any prior knowledge. If this initial ensemble is too far from the truth, the Euler forward integration can fail in this spin-up period. This fact is illustrated with a simple experiment illustrated in figure 4. The (unstable) steady state of L96 is $\mathbf{X}_j = F \ \forall j$, where $F$ is the forcing term in (33). A simple way to generate an initial ensemble (in the absence of any prior information) is to add random perturbations to this steady state for each one of the $N = 40$ variables and $M = 10$ ensemble members. Using DETKBF with $\delta = 0.05$ and $\lambda = 4$, we generated initial ensembles using the multiples of the observational error covariance: $\mathbf{R}$ (left panel), $2\mathbf{R}$ (center panel) and $3\mathbf{R}$ (right panel); we show the first 150 analysis cycles. The Euler forward method can lead to large initial increase in analysis RMSE before the filter stabilizes. This does not happen with the DSI method. Moreover, for larger initial perturbations the Euler forward fails (not shown).

### 5.3. The SPEEDY model

Finally, we implement ETKBF and DETKBF in a model that is more representative of those used in operational numerical weather prediction. We choose a medium-complexity AGCM developed by Molteni (2003) and known as SPEEDY (Simplified Parameterizations, primitivE-Equation Dynamics) which has been adapted for data assimilation by Miyoshi (2005). This model has a spectral primitive-equation dynamic core and a set of simplified physical parameterization schemes; it achieves computational efficiency while maintaining realistic simulations similar to those of state-of-the-art AGCMs with complex physics. The model has a resolution of T30L7, with horizontal spectral truncation at 30 wave numbers and 7 vertical levels. Data are output on a horizontal grid of 96 longitudinal and 48 latitudinal points. The SPEEDY model is formulated in $\sigma$-coordinates and calculates five field variables: zonal wind $u$, meridional wind $v$, temperature $T$, relative humidity $q$ and surface pressure $ps$. The geopotential height $z$ for different pressure levels may be obtained by interpolation.

The nature run starts after a one-year spin-up from state of rest. The "observations" are generated by adding Gaussian random perturbations to every variable with the following standard deviations: $1 m/s$ for $u$ and $v$, $1 K$ for $T$, $10^{-3} kg_{water}/kg_{air}$ for $q$ and $1 hPa$ for $ps$. Observations are taken every 6 hr at all 7 vertical levels at horizontal positions that resemble a realistic radiosonde observational network (figure 5). The density of observations is larger over continents than over the oceans (which are barely observed), and the Northern Hemisphere is better observed than the Southern Hemisphere.

For the assimilation, an ensemble of $M = 20$ members is used. The R-localization parameters are $\lambda = 500\,km$ in the horizontal and $\lambda_v = 0.1\ln p$ in the vertical. As with L96, adaptive multiplicative covariance inflation (Miyoshi, 2011) is applied to avoid manual tuning. For ETKBF and DETKBF we look for the minimum number of pseudo-time steps that lead to a performance equal to that of the LETKF. Using 3 steps or less leads to noticeable differences, starting at 4 the differences are minimal and by 6 the impact is practically indistinguishable. For this number of steps the computational time required for an assimilation cycle is comparable to that of LETKF.

To illustrate the effectiveness of ETKBF and DETKBF, figure 6 shows the results of a single-observation experiment for observations at two locations: one over the Labrador Peninsula in a well-observed region, and one over the Southern Pacific in a poorly-observed region (figure 5). These locations are used for single-observation experiments to illustrate the equivalence in the effect of the LETKF and the Bucy-type formulations with 6 steps; the effect is the same when all the observations are used. The first row of figure 6 shows the result of the single-observation experiment using the observation over the Labrador Peninsula. We illustrate the result of the assimilation as the difference between analysis mean and background mean (i.e. $\bar{\mathbf{x}}^a - \bar{\mathbf{x}}^b$) for the zonal wind at 510 hPa. This update obtained by the three filters presents no differences. The second row shows the result of the experiment using the observation over the Southern Pacific using the ratio of analysis ensemble spread over background ensemble spread. Since this is a poorly-observed region, the reduction is considerable, e.g. there are regions where this ratio is as low as 0.3. This large reduction implies large values of $\beta$ (section 3). In spite of this condition, the results obtained from the three filters are the same.

To assess the performance of the three assimilation techniques, a latitude-weighted RMSE is computed for each one of the variables:

$$RMSE = \sqrt{\frac{1}{J}\sum_{j=1}^{J}(x^a_i - x^t_i)^2 \cos\varphi_i} \qquad (32)$$

where $x^t$ corresponds to the nature run, $x^a$ corresponds to the analysis, and $\varphi$ is the latitude angle. We consider each one of the 6 variables ($u$, $v$, $T$, $q$, $ps$, $z$) separately at each one of the 7 vertical levels. In figure 7 we show the analysis RMSE for temperature at three different vertical levels: 950hPa in the left column, 510 hPa in the center column and 200hPa in the right column. The RMSE is first computed globally (top row), and then by region: Northern Hemisphere (second row from top), tropics (third row from top) and Southern Hemisphere (bottom row). In each panel the bars represent 1 standard deviation around the mean. The performance of the three filters is indistinguishable; the computing time was comparable as well. The effect of the observational density is clear (and the same for the three filters): e.g. for the (well-observed) NH the mean analysis RMSE is about half of the observational error while for the SH it is above this value. The results are similar for all the variables at all sigma levels (not shown). The DSI method ensures a proper performance even in the data sparse regions where the stiffening of the ODEs may happen.

## 6. Summary and conclusions

We have analyzed two recently proposed ensemble formulations based on the Kalman-Bucy filter. In BGR09, the ensemble of perturbations is updated by the solution of an ODE in pseudo-time and the mean is updated with the standard Kalman equations. In BR10 the full ensemble is updated by a solution of a single set of ODEs.

In this work we have achieved two objectives. First, a study of the ODEs involved in these Bucy-type formulations was performed. The ratio of the magnitude of the background covariance with respect to the magnitude of the observational error covariance is crucial for the behavior of these ODEs. The ODEs stiffen under certain conditions and cause the failure of the Euler forward integration used in these works. As an alternative, a diagonal semi-implicit integration method with variable step size was introduced; this method ensures stability and is computationally affordable.

Second, transform-based versions of BGR09 and BR10 were developed. For these alternatives, the variables integrated in pseudo-time are weights, with dimension equal to the ensemble size rather than the much larger model dimension. The availability of the weights is useful for some applications (Kalnay and Yang, 2010, Yang and Kalnay, 2009a).

Three models were used to test our transform-based Kalman-Bucy filters and to illustrate the value of the diagonal eemi-implicit integration method. First, the L63 model allowed us to perform experiments with frequent and infrequent observations. For frequent observations, it was found that with 3 steps the Bucy-type formulations achieve performances comparable to that of the LETKF; with 5 steps the performances were indistinguishable. For infrequent observations, the proposed Diagonal Semi-Implicit method with variable time stepping proved its value: to achieve the performance of the LETKF only 8 pseudo-time steps were needed, a huge reduction from the number needed using EF (50 for ETKBF and 300 for DETKBF).

In the L96 model we applied ETKBF and DETKBF using R-localization and adaptive multiplicative covariance inflation. It was found that the performance of these schemes was equal to that of LETKF with only 4 pseudo-time steps. The advantages of using the DSI integration when initializing the background ensemble without prior information were demonstrated with an example.

Finally, we implemented our schemes in an AGCM known as the SPEEDY model with a realistic radiosonde observational network. The equivalence in their performance with respect to the LETKF was shown, even for data sparse regions (e.g. over the oceans) in which the ODEs are bound to stiffen.

An essential implementation issue for the Bucy-type formulations is the choice of the number of steps for the integration. We have shown that in the 'frequent observations' case (corresponding to $\Delta t = 0.08$ in L63, $\Delta t = 0.05$ in L96, and $\Delta t = 6 hr$ in an AGCM), an adequate performance starts at 3-6 steps. For infrequent observations ($\Delta t = 0.25$ in L63) this number doesn't surpass 10 (as a result of using the DSI method). For any dynamical system, it will be necessary to first estimate $\beta$ as in (22) for the given assimilation window length. A possible improvement of the R-localization implementation would to compute $\beta$ locally and let every gridpoint use a different number of steps depending on the local degree of stiffness.

The computational implementation of the Bucy-type approaches and their transform versions are straightforward and amenable to parallel computing. Finally, the continuous formulation of the ensemble Kalman filter allows for a seamless implementation of the incremental analysis update (IAU, Bloom et al, 1996) as demonstrated in the mollified ensemble Kalman-Bucy filter (Bergemann and Reich, 2010a). The purpose of this implementation is to avoid the imbalance introduced by the jumps from background to analysis that are present in sequential data assimilation; its performance will be tested in a forthcoming study.

## 7. Acknowledgements

The authors are thankful for the valuable input from Dr. John Thuburn and two anonymous reviewers; their comments improved the format and content of this work. The support of NASA grant

NNX07AM97G, NNX08AD40G, DOE grant DEFG0207ER64437, NOAA grant NA09OAR4310178, and ONR grant N000140910418, N000141010557 is gratefully acknowledged.

**Appendix A: Derivation of BR10 and DETKBF**

The ODE representing the analysis step for the full ensemble can be written as:

$$\frac{d\bar{\bar{\mathbf{X}}}}{ds} = \frac{d}{ds}\left(\mathbf{X} + \bar{\mathbf{x}}\mathbf{1}^T\right) = \frac{d}{ds}\mathbf{X} + \left(\frac{d}{ds}\bar{\mathbf{x}}\right)\mathbf{1}^T \qquad (A1)$$

The two terms of (A1) correspond to the ensemble version of the KBF equations for the perturbations and the mean, respectively, i.e.,

$$\frac{d\mathbf{X}}{ds} = -\frac{1}{2(M-1)}\mathbf{X}\mathbf{X}^T\mathbf{H}^T\mathbf{R}^{-1}\mathbf{H}\mathbf{X} \qquad (A2)$$

$$\frac{d\bar{\mathbf{x}}}{ds} = -\frac{1}{M-1}\mathbf{X}\mathbf{X}^T\mathbf{H}^T\mathbf{R}^{-1}(\mathbf{H}\bar{\mathbf{x}} - \mathbf{y}) \qquad (A3)$$

Substitute (A2) and (A3) into (A1) and factorize:

$$\frac{d\bar{\bar{\mathbf{X}}}}{ds} = -\frac{1}{M-1}\mathbf{X}\mathbf{X}^T\mathbf{H}^T\mathbf{R}^{-1}\left[\frac{1}{2}\mathbf{H}\mathbf{X} + (\mathbf{H}\bar{\mathbf{x}} - \mathbf{y})\mathbf{1}^T\right] \qquad (A4)$$

(A4) describes the assimilation of the $\bar{\bar{\mathbf{X}}}$ in terms of $\mathbf{X}$ and $\bar{\mathbf{x}}$. To obtain an expression in terms of $\bar{\bar{\mathbf{X}}}$, recall that $\bar{\mathbf{x}} = M^{-1}\bar{\bar{\mathbf{X}}}\mathbf{1}$ and $\mathbf{X} = \bar{\bar{\mathbf{X}}}(\mathbf{I} - \mathbf{U})$, where $\mathbf{U} = M^{-1}\mathbf{1}\mathbf{1}^T$. Substituting these into (A4), using the fact that $\mathbf{I} - \mathbf{U}$ is symmetric and idempotent and simplifying we obtain equation (15):

$$\frac{d\bar{\bar{\mathbf{X}}}}{ds} = -\frac{1}{M-1}\bar{\bar{\mathbf{X}}}(\mathbf{I} - \mathbf{U})\bar{\bar{\mathbf{X}}}^T\mathbf{H}^T\mathbf{R}^{-1}\left[\frac{1}{2}\mathbf{H}\bar{\bar{\mathbf{X}}}(\mathbf{I} + \mathbf{U}) - \mathbf{y}\mathbf{1}^T\right] \qquad (A5)$$

**Appendix B. Equivalence of ETKBF and LETKF for an infinite number of steps in pseudo-time**

It can be shown that in the limit of infinite steps, the ETKBF is numerically equivalent to the LETKF (Hunt *et al.*, 2007). To proceed, we start by writing the expressing the covariance in ensemble space for any instant in pseudo-time:

$$\tilde{\mathbf{P}}(s) = \frac{\mathbf{W}(s)\mathbf{W}^T(s)}{M-1} \tag{B1}$$

with $\tilde{\mathbf{P}}(0) = \tilde{\mathbf{P}}^b = (M-1)^{-1}\mathbf{I}$ and $\tilde{\mathbf{P}}(0) = \tilde{\mathbf{P}}^a$. Using the chain rule, we can find the pseudo-time derivate of this expression and perform simplifications:

$$\frac{d}{ds}\tilde{\mathbf{P}}(s) = \frac{2}{M-1}\left(\frac{d}{ds}\mathbf{W}(s)\right)\mathbf{W}^T(s) = \frac{2}{M-1}\left(-\frac{1}{2(M-1)}\mathbf{W}(s)\mathbf{W}^T(s)\mathbf{Y}^{bT}\mathbf{R}^{-1}\mathbf{Y}^b\mathbf{W}(s)\right)\mathbf{W}^T(s)$$

Hence,

$$\frac{d}{ds}\tilde{\mathbf{P}} = -\frac{\mathbf{W}\mathbf{W}^T}{M-1}\mathbf{Y}^{bT}\mathbf{R}^{-1}\mathbf{Y}^b\frac{\mathbf{W}\mathbf{W}^T}{M-1} = -\tilde{\mathbf{P}}\mathbf{Y}^{bT}\mathbf{R}^{-1}\mathbf{Y}^b\tilde{\mathbf{P}} \tag{B2}$$

The analytical solution to this Riccati equation is (see e.g. Simon 2006):

$$\tilde{\mathbf{P}}(s) = \tilde{\mathbf{P}}(0)\left(\mathbf{Y}^{bT}\mathbf{R}^{-1}\mathbf{Y}^b\tilde{\mathbf{P}}(0)s + \mathbf{I}\right)^{-1} = \frac{\mathbf{I}}{M-1}\left(\mathbf{Y}^{bT}\mathbf{R}^{-1}\mathbf{Y}^b\frac{s}{M-1} + \mathbf{I}\right)^{-1}$$

$$\tilde{\mathbf{P}}(s) = \left(\mathbf{Y}^{bT}\mathbf{R}^{-1}\mathbf{Y}^b s + (M-1)\mathbf{I}\right)^{-1} \tag{B3}$$

In particular, for $s = 1$ we get the same expression as for the LETKF:

$$\tilde{\mathbf{P}}(1) = \tilde{\mathbf{P}}^a = \left(\mathbf{Y}^{bT}\mathbf{R}^{-1}\mathbf{Y}^b + (M-1)\mathbf{I}\right)^{-1} \tag{B4}$$

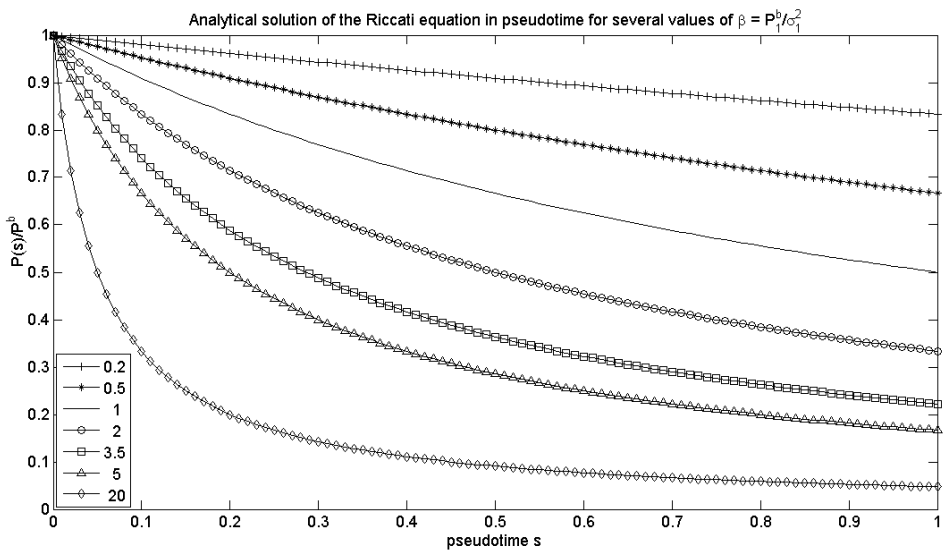

Figure 1. Analytical solution of the Bucy covariance equation in pseudo-time for a scalar case with the variable observed directly. Different lines correspond to different values ratio of background variance over observational error variance (as quantified in the legend). The ODE stiffens as this ratio becomes larger.

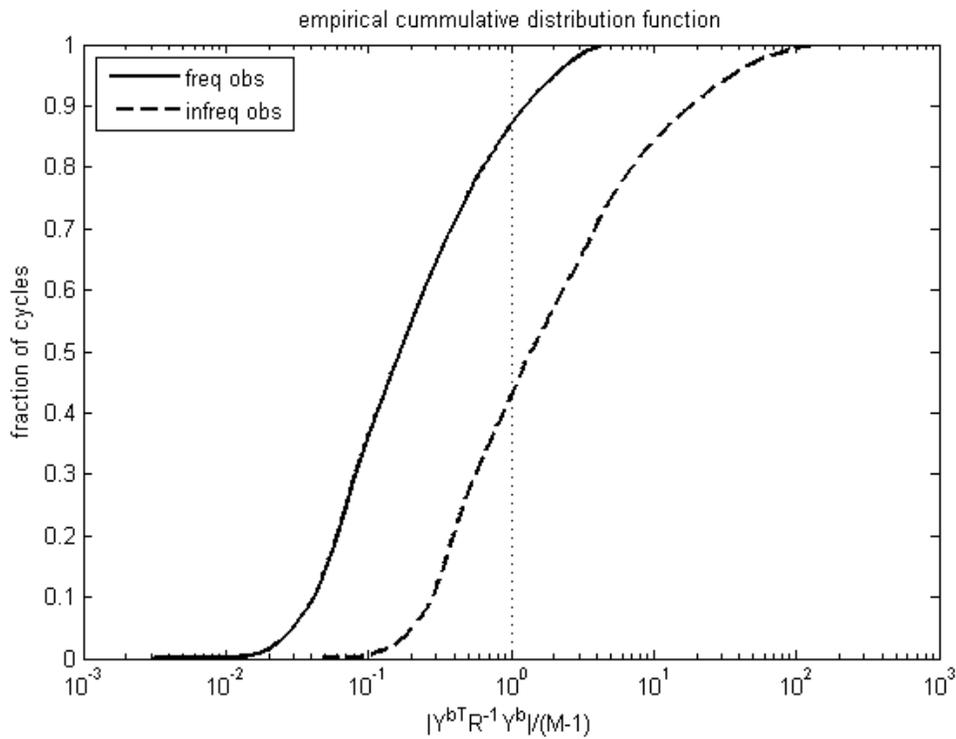

Figure 2. Empirical cumulative distribution function of $\beta = \dfrac{\left|\mathbf{Y}^{b^T}\mathbf{R}^{-1}\mathbf{Y}^b\right|}{M-1}$ for short and long assimilation windows using the Lorenz 1963 model. The value of this ratio for infrequent observations is in general an order of magnitude larger than for frequent observations (details explained in section 5.1).

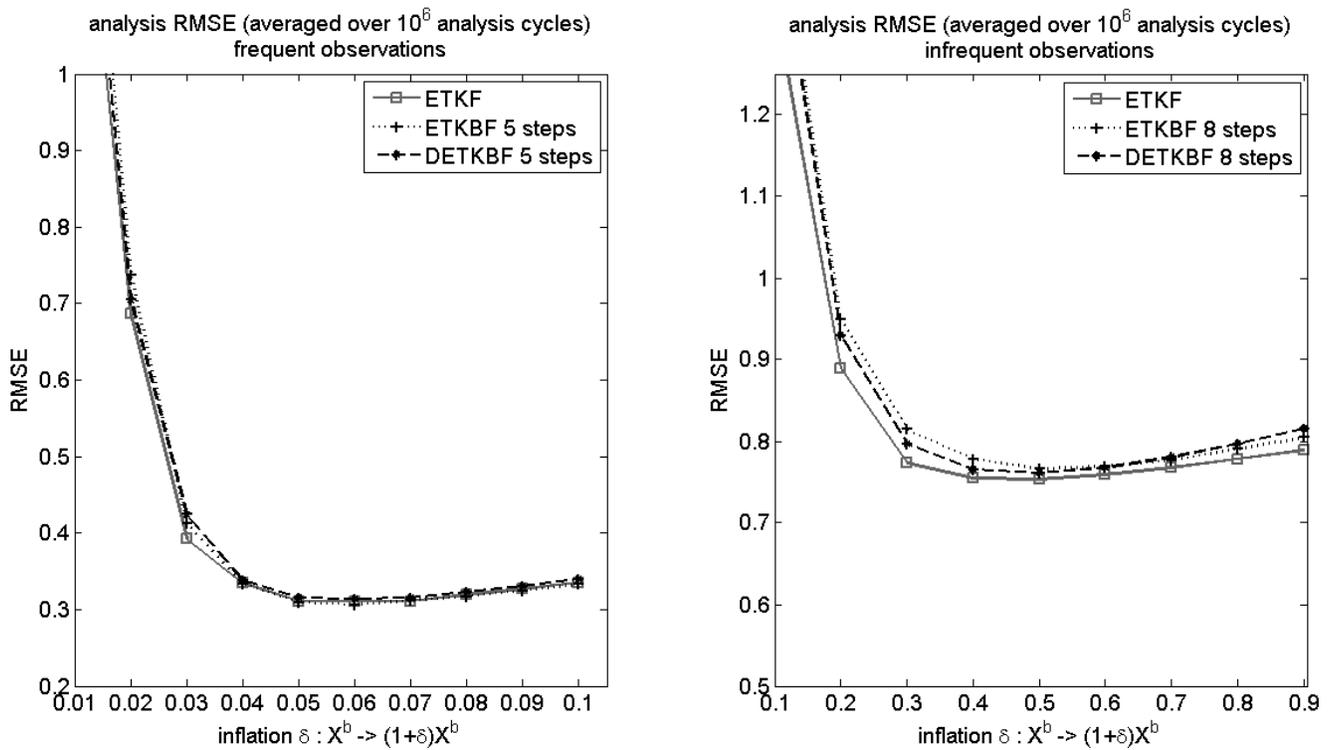

Figure 3. Analysis RMSE values (averaged over $10^6$ assimilation cycles) for ETKF, ETKBF and DETKBF in L63, as a function of multiplicative inflation. In the left panel, the frequent observations case is shown; the Bucy-type formulations were integrated using the diagonal semi-implicit (DSI) method with 5 uniform steps. In the right panel, the infrequent observations case is shown; the Bucy-type formulations were integrated using the (DSI) method with 8 non-uniform steps.

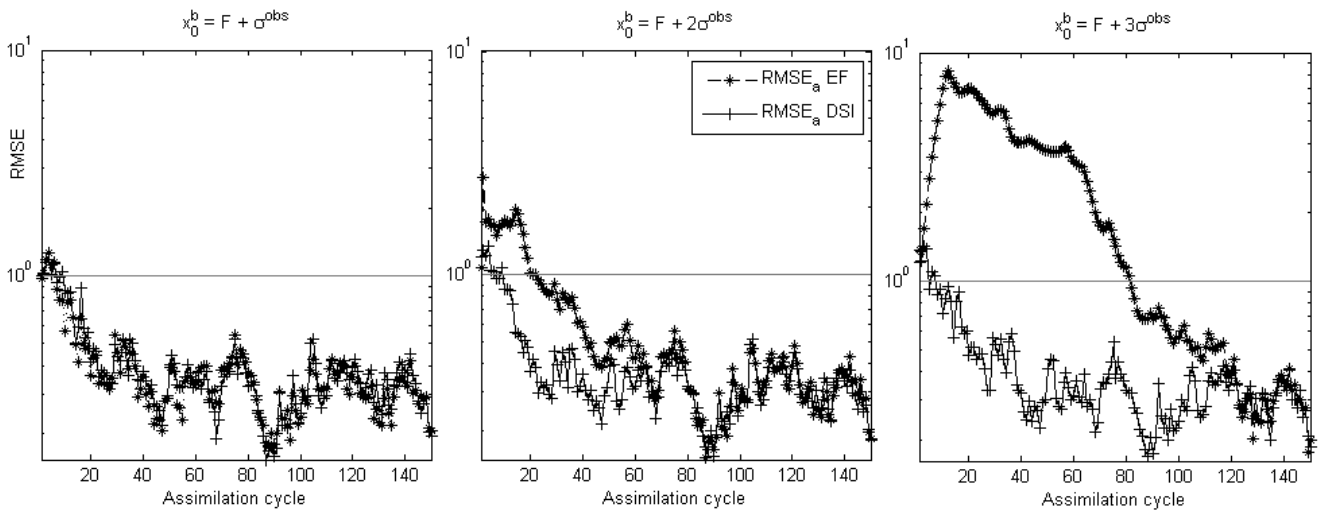

Figure 4. Analysis RMSE for the first 150 assimilation cycles of an experiment using L96. The effect of the two integration schemes (Euler forward and diagonal semi-implicit) is shown for different initial ensembles. As the initial ensemble is more inaccurate (from left to right), EF takes longer to initialize the filter while DSI does not present problems.

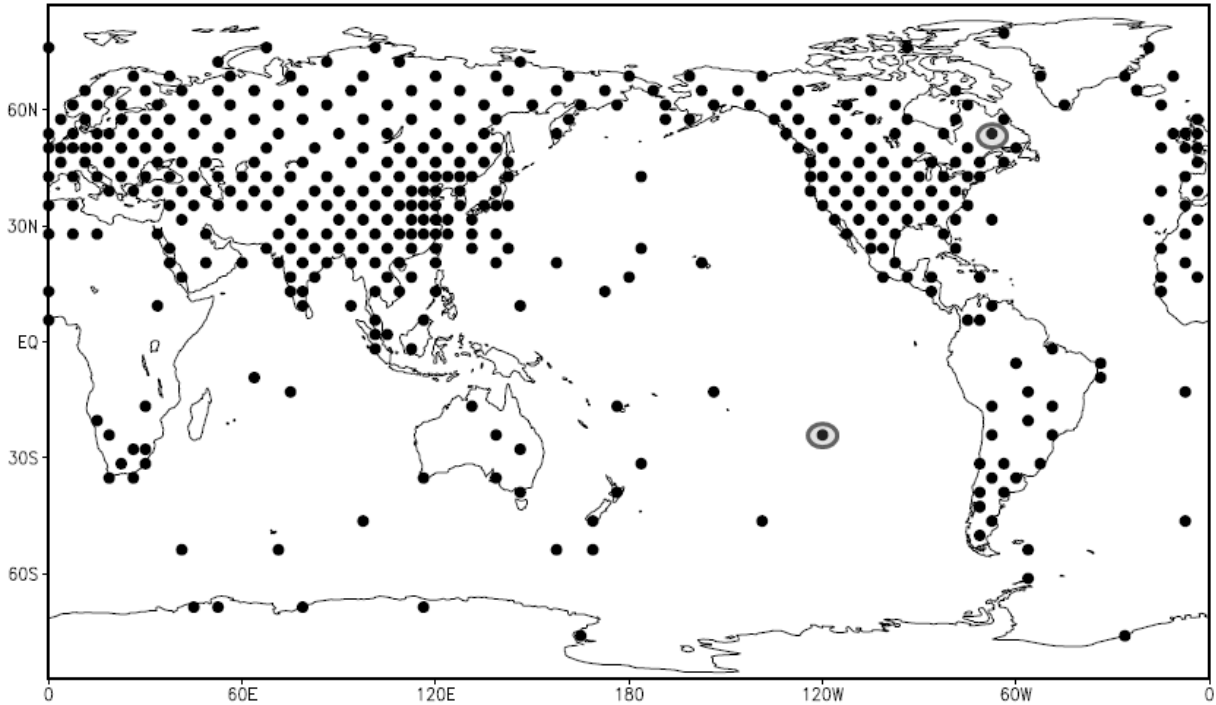

Figure 5. Observational density used for the experiments with the SPEEDY model. The spatial distribution of the observations resembles a realistic radiosonde network (Miyoshi, 2011). Two positions are circled, one over the Labrador Peninsula and other over the Southern Pacific Ocean; these are used for the experiments depicted in figure 6.

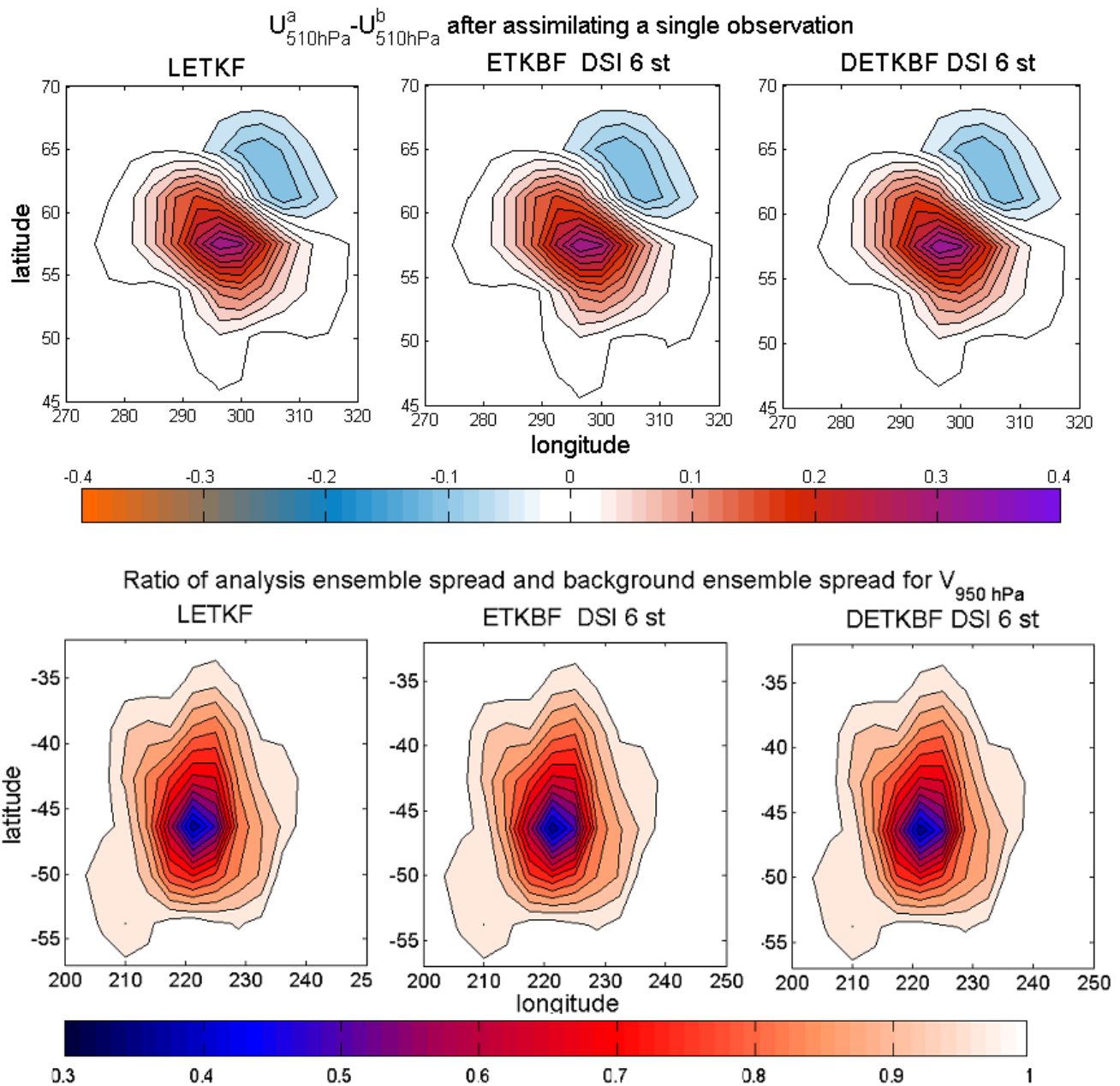

Figure 6. Impact of two single-observation experiments. In the first row, the observation is located in the well-observed Labrador Peninusula. The difference between analysis mean and background mean is the same using the three methods; the variable illustrated is zonal wind at 510 hPa. In the second row, the observation is located in the poorly-observed Southern Pacific. In this case we show the ratio of the analysis spread to background spread for the meridional wind at 950 hPa; the performance of the three filters is indistinguishable.

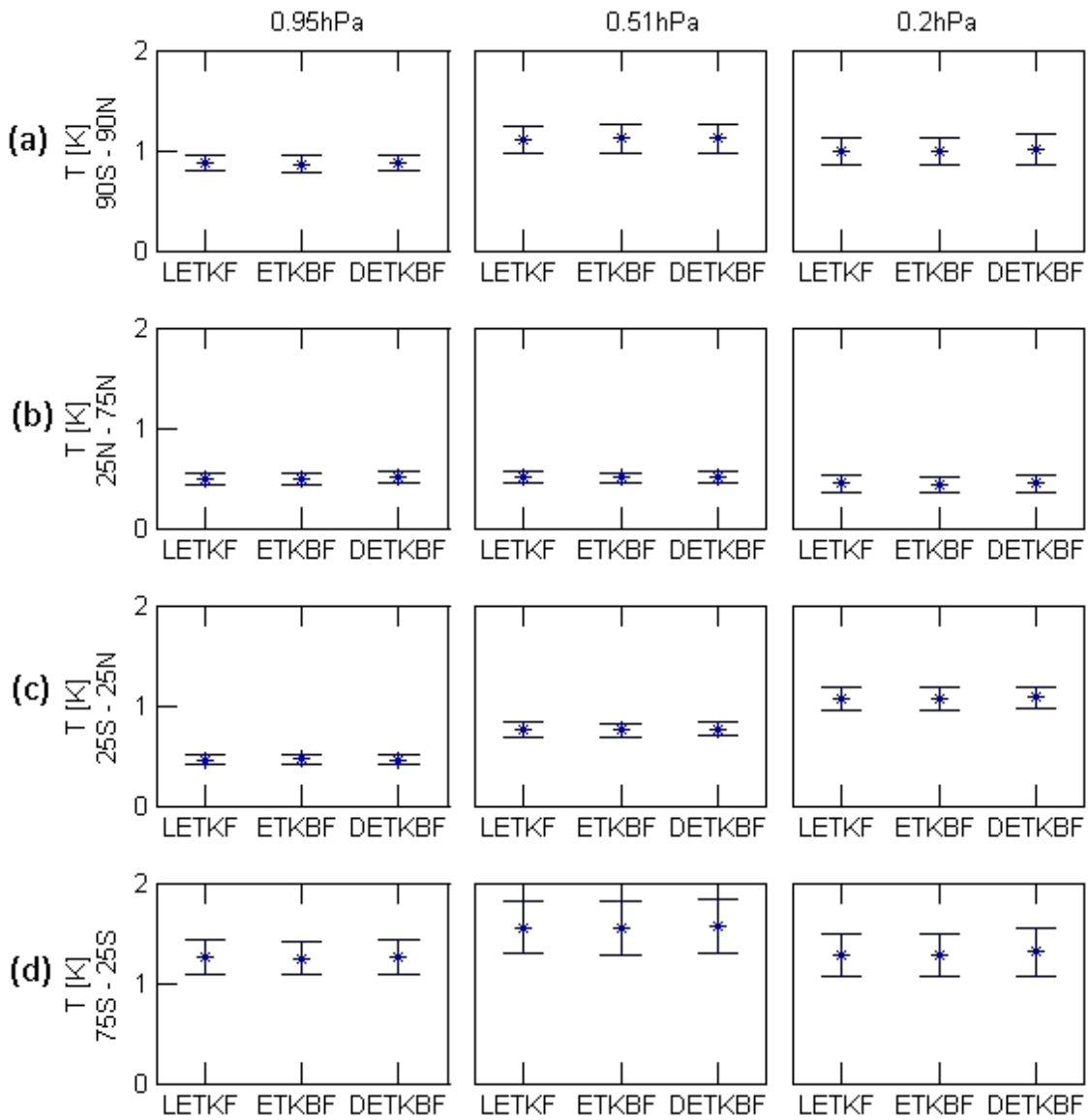

Figure 7. Analysis RMSE for the temperature in January using LETKF, ETKBF and DETKBF. Each column shows a different vertical level. This diagnostic was computed (a) globally, (b) for the Northern Hemisphere, (c) for the tropics and (d) for the Southern Hemisphere. The performance of the three filters is indistinguishable (the bars represent 1 standard deviation around the mean).

|  | LETKF | ETKBF (4 steps) | DETKBF (4 steps) |
|---|---|---|---|
| $M = 10$ | RMSE = 0.3215 (0.0832) <br> spread = 0.3532(0.0324) <br> $\delta$ = 0.0289 (0.0112) | RMSE = 0.3215 (0.0862) <br> spread = 0.3515(0.0327) <br> $\delta$ = 0.0289 (0.0114) | RMSE = 0.3227 (0.0883) <br> spread = 0.3513(0.0330) <br> $\delta$ = 0.0289 (0.0115) |
| $M = 15$ | RMSE = 0.3190 (0.0789) <br> spread = 0.3694(0.0329) <br> $\delta$ = 0.0294 (0.0115) | RMSE = 0.3184 (0.0793) <br> spread = 0.3671(0.0326) <br> $\delta$ = 0.0292 (0.0114) | RMSE = 0.3197 (0.0791) <br> spread = 0.3670(0.0328) <br> $\delta$ = 0.0292 (0.0114) |

Table 1. Results of the experiments with the L96 model. Three assimilation methods (columns) and two ensemble sizes (rows) are used. In each case, a sample of $10^6$ assimilation cycle was used to find the average analysis RMSE, average ensemble spread and average inflation parameter (averaged also over the 40 gridpoints). The numbers in parenthesis correspond to the standard deviations of the reported parameters.